\documentclass{article}
\usepackage{amsmath}
\usepackage{amsfonts}
\usepackage{amssymb}	
\usepackage{color}
\newtheorem{theo}{Theorem}[section]
\newtheorem{lem}[theo]{Lemma}
\newtheorem{cor}[theo]{Corollary}
\newtheorem{rem}[theo]{Remark}

\newcommand{\mysection}[1]{\section{#1} \setcounter{equation}{0}}
\newcommand{\proof}{{\sc Proof.} \quad}

\newcommand{\R}{\mathbb{R}}

\newcommand{\be}{\begin{equation} \label}
\newcommand{\ee}{\end{equation}}
\newcommand{\bes}{\begin{equation} \begin{array}{c} \label}
\newcommand{\ees}{\end{array} \end{equation}}
\newcommand{\bea}{\begin{eqnarray}\label}
\newcommand{\eea}{\end{eqnarray}}
\newcommand{\beas}{\begin{eqnarray} \begin{array}{rcl} \label}
\newcommand{\eeas}{\end{array} \end{eqnarray}}
\newcommand{\bas}{\begin{eqnarray*}}\newcommand{\eas}{\end{eqnarray*}}
\newcommand{\bass}{\begin{eqnarray*} \begin{array}{rcl}}
\newcommand{\eass}{\end{array} \end{eqnarray*}}
\newcommand{\basss}{\begin{eqnarray*} \begin{array}{c}}
\newcommand{\easss}{\end{array} \end{eqnarray*}}
\newcommand{\qed}{{}\hfill $\square$ \\}
\newcommand{\bit}{\begin{itemize}}
\newcommand{\eit}{\end{itemize}}
\newcommand{\nn}{\nonumber}
\newcommand{\eps}{\varepsilon}
\newcommand{\abs}{\\[3mm]}

\newcommand{\pO}{\partial\Omega}

\newcommand{\io}{\int_\Omega}

\newcommand{\F}{{\cal F}}
\newcommand{\D}{{\cal D}}
\newcommand{\B}{{\cal B}}
\newcommand{\set}{{\cal S}(m,M,B)}
\newcommand{\tm}{T_{max}(u_0,v_0)}

\begin{document}
\title{Finite-time blowup in a supercritical quasilinear parabolic-parabolic Keller-Segel system in dimension $2$}
\author{
Tomasz Cie\'{s}lak \\
{\small Institute of Mathematics, Polish Academy of Sciences, \'Sniadeckich 8, 00-956 Warsaw, Poland}\\
{\small E-Mail: T.Cieslak@impan.pl}
\and 
Christian  Stinner\\
{\small Institut f\"ur Mathematik, Universit\"at Z\"{u}rich, Winterthurerstrasse 190, 8057 Z\"{u}rich, Switzerland} \\
{\small E-Mail: christian.stinner@math.uzh.ch}
}
\date{}
\maketitle
\begin{abstract}
\noindent
  In this paper we prove finite-time blowup of radially symmetric solutions to the quasilinear parabolic-parabolic two-dimensional Keller-Segel system for any positive mass. This is done in case of nonlinear diffusion and also in the case of nonlinear cross-diffusion provided the nonlinear chemosensitivity term is assumed not to decay. Moreover, it is shown that
the above-mentioned lack of non-decay assumption is essential with respect to keeping the dichotomy finite-time blowup
against boundedness of solutions. Namely, we prove that without the non-decay assumption possible asymptotic behaviour of solutions includes also infinite-time blowup. 

\end{abstract}
\mysection{Introduction}\label{section1}
In the present paper we deal with solutions $(u,v)$ of the parabolic-parabolic 
Keller-Segel system
\be{0}
	\left\{ \begin{array}{ll}
	u_t= \nabla \cdot (\phi(u) \nabla u) - \nabla \cdot ( \psi(u) \nabla v ), & \; x\in\Omega, \ t>0, \\[2mm]
	v_t=\Delta v-v+u, & \; x\in\Omega, \ t>0, \\[2mm]
	\frac{\partial u}{\partial\nu}=\frac{\partial v}{\partial\nu}=0, & \; x\in\partial\Omega, \ t>0, \\[2mm]
	u(x,0)=u_0(x), \quad v(x,0)=v_0(x), & \; x\in\Omega,
	\end{array} \right.
\ee
in a ball $\Omega = B_R \subset \R^2$, $R>0$,
where the
initial data are supposed to satisfy $u_0 \in C^0(\bar\Omega)$ and $v_0\in W^{1,\infty}(\Omega)$ such that
$u_0 > 0$ and $v_0 > 0$ in $\bar{\Omega}$.

Moreover, let $\phi, \psi \in C^2([0,\infty))$ such that  
\begin{equation}\label{0.1}
 \phi(s) >0, \qquad \psi(s) = s \beta (s), \quad\mbox{and }\quad \beta (s) >0 \quad\mbox{for } s \in [0,\infty) 
\end{equation}
are fulfilled with some $\beta \in C^2 ([0,\infty))$.

Let us introduce the following notation. Suppose that there exist $s_0 >1$ and positive constants $a$ and $b$ such that the
functions
\begin{equation}\label{GH}
 G(s) := \int\limits_{s_0}^s \int\limits_{s_0}^\sigma \frac{\phi(\tau)}{\psi(\tau)} \; d\tau \, d \sigma, \quad s>0,
 \qquad\mbox{and}\qquad 
 H(s) := \int\limits_0^s \frac{\sigma \phi(\sigma)}{\psi(\sigma)} \; d\sigma, \quad s \ge 0,
\end{equation}
satisfy
\begin{equation}\label{G1}
 G(s) \le a s(\ln s)^\mu, \; s \ge s_0,  
\end{equation}
with some $\mu \in (0,1)$ as well as
\begin{equation}\label{H1}
  H(s) \le b \frac{s}{\ln s}, \; s \ge s_0.
\end{equation}
We remark that $H$ in \eqref{GH} is well-defined due to the positivity of $\beta$ in $[0,\infty)$.

Moreover, assume that 
\begin{equation}\label{psi}
  \psi(s) \ge c_0 \, s, \quad s \ge 0.
\end{equation}
Next we introduce the well-known Liapunov functional for the Keller-Segel system. 
\begin{equation}\label{F}
  \F(u,v):=\frac{1}{2} \io |\nabla v|^2 + \frac{1}{2} \io v^2 - \io uv + \io G(u)
\end{equation}
is a Liapunov functional for \eqref{0} with dissipation rate
\begin{equation}\label{D}
  \D(u,v):=\io v_t^2 + \io \psi(u) \cdot \Big| \frac{\phi (u)}{\psi(u)} \nabla u - \nabla v \Big|^2.
\end{equation}
More precisely, any classical solution to \eqref{0} satisfies
\begin{equation}\label{liapunov}
	\frac{d}{dt} \F(u(\cdot,t),v(\cdot,t)) = - \D (u(\cdot,t),v(\cdot,t)) 
	\qquad \mbox{for all } t \in (0,\tm),
\end{equation}
where $\tm \in (0,\infty]$ denotes the maximal existence time of $(u,v)$ (see \cite[Lemma~2.1]{win_mmas}). \abs

Our main result is a finite-time blowup in the case of a quasilinear problem provided the nonlinear chemosensitivity term 
satisfies $\psi(u)\geq Cu^q$ with some $q\geq 1$ and $C>0$. The considered system \eqref{0} was introduced in \cite{KS:1} to describe the motion of cells on a surface, where the cells are diffusing and moving towards the gradient of a substance called chemoattractant, the latter being produced by the cells themselves. The main motivation was to describe the chemotactic collapse of cells interpreted as finite-time blowup of the component $u$ of a solution to \eqref{0}. However, almost all results concerning the finite-time blowup of solutions to \eqref{0} were proved for its parabolic-elliptic simplification. Main achievements concerning this issue are contained in \cite{jl:expl,B,Na95} stating the chemotactic collapse for a semilinear system, i.e. $\phi(u)=1$ and  $\psi(u)=u$, provided that the initial mass exceeds $8\pi$ in the radially symmetric case or $4\pi$ in the case of solutions without the assumption of symmetry, see \cite{Nagai}.  Moreover, it has been shown that in higher dimensions a finite-time blowup of solutions to the semilinear version of \eqref{0} can occur independently of the initial mass provided that the initial data are concentrated enough \cite{Na95}.  Finally, in the case of a quasilinear system, for any space dimension $n$ critical nonlinearities have been identified such that if $\phi$ and $\psi$ satisfy the subcritical relation, then solutions to \eqref{0} stay bounded for any time, while for those satisfying the supercritical relation solutions blow up in finite-time independently of the magnitude of initial mass provided the data are concentrated enough, see \cite{djiewin}.

However, all those results are available only for a parabolic-elliptic simplification of \eqref{0}. In the case of the original fully parabolic version the investigation of chemotactic collapse turned out to be a much more challenging issue. So far the only existing result in the literature showing the occurrence of finite-time blowup of solutions to the semilinear version of \eqref{0} in space dimension $2$ is the one in \cite{HV:1}, where an example of a special solution to the semilinear version of \eqref{0} in dimension $n=2$ blowing up in a finite time is shown. Moreover, there are several results concerning other dimensions. In \cite{clKS} the explosion of solutions to the one-dimensional Keller-Segel system with appropriately weak diffusion of cells, properly large mass and sufficiently fast diffusion of chemoattractant is shown. In \cite{win_bu} M. Winkler introduced a new method which led him to the finite-time blowup of solutions to the semilinear Keller-Segel system in dimensions $n\geq 3$. His method was generalized in \cite{ciesti} and the result was extended to the quasilinear case for the optimal range of nonlinearities.
This way, to the best of our knowledge, we present a first result concerning a finite-time blowup of solutions to the fully parabolic quasilinear Keller-Segel system in dimension two. Moreover, we show that as expected an explosion 
takes place independently of the size of initial mass. This result is proved in both cases $\psi(u)=u$ (nonlinear 
diffusion case) and a fully nonlinear cross-diffusion.  Both results are optimal in view of possible nonlinearities generating finite-time blowup. As the result in \cite{ss} shows we prove finite-time blowup for the optimal range of nonlinear diffusion in the case $\psi(u)=u$. Moreover, at least under the restriction of polynomial nonlinearities we have the optimal result provided we accept the non-decay assumption on nonlinear chemotactic sensitivity, see \cite{taowin_jde}.
On the other hand, again assuming the non-decay of $\psi$, even without assuming nonlinearities to be polynomial we still have the exhaustive finite-time blowup result, see \cite{ja2}. 
Next, it is shown that the above-mentioned non-decay assumption is essential for finding critical exponents distinguishing between finite-time blowup and boundedness in the case of nonlinear cross-diffusion.  Despite the fact that nonlinearities considered in that case seem to be from the finite-time blowup regime, we construct solutions to \eqref{0} in dimension $2$ which blow up in infinite time when the nonlinear chemotactic sensitivity term is decreasing.  

For any $\phi$ and $\psi$ from the class defined in the beginning we have the local existence of smooth solutions. 
Moreover, the solution $(u,v)$ is positive for $t>0$ and preserves mass as well as radial symmetry. In
particular, it satisfies 
\begin{equation}\label{Mduze}
\int_\Omega u(x,t) dx = \int_\Omega u_0 (x) dx \qquad \mbox{and}\qquad
\int_\Omega v(x,t)dx\leq \max\left\{\int_\Omega u_0(x)dx,\int_\Omega v_0(x)dx\right\}
\end{equation} 
for all $t\in(0,T_{max}(u_0,v_0))$, where $T_{max}(u_0,v_0)$ is the maximal time of existence of solutions. Furthermore, for solutions which cease to exist
for all positive times, $\| u (\cdot, t) \|_{L^\infty (\Omega)}$ blows up in finite time. 
For the details we refer to \cite[Lemma 2.1]{ciesti}. 

Our main theorems are the following.
\begin{theo}\label{theo1}
  Suppose that $\Omega=B_R \subset \R^2$ with some $R>0$, assume also that \eqref{G1}, \eqref{H1}, and \eqref{psi}
  are satisfied. Next let $m>0$ and $A>0$ be given. 
  Then there exist positive constants $T(m,A)$ and $K(m)$ such that for any
  \begin{eqnarray}\label{t1.1}
	(u_0,v_0) \in \B (m,A) &:=& \bigg\{
	(u_0,v_0) \in C^0(\bar\Omega) \times W^{1,\infty}(\Omega) \ \bigg| \
	\mbox{$u_0$ and $v_0$ are radially symmetric} \nn \\ 
    & & \hspace*{5mm} \mbox{and positive in $\bar\Omega$, 
	$\io u_0=m$, $\|v_0\|_{W^{1,2}(\Omega)} \le A$,}\nn \\ 
    & & \hspace*{5mm} \mbox{and $\F(u_0,v_0) \le -K(m) \cdot (1+A^2)$} \bigg\},
  \end{eqnarray}
  the corresponding solution $(u,v)$ of \eqref{0} blows up at the finite time $\tm \in (0,\infty)$, where
  $\tm \le T(m,A)$. Furthermore, for any $m>0$ there exists $A>0$ such that the set $\B (m,A)$ is nonempty.
\end{theo}
Next let us introduce the following corollary simplifying our result in the case of $\psi(u):=u$. It covers the interesting case of a system with nonlinear diffusion. The result is optimal in view of its global existence counterpart proved in \cite{ss}.

\begin{cor}\label{cor1.3}
  Assume that $\psi(s) =s$ for $s \ge 0$ and that  
  $\phi(s)\leq Cs^{q}$, $s \ge 1$, for some $q<0$ and $C>0$. Furthermore, suppose that  
  $\phi$ is a decreasing function. Let
  $\Omega=B_R \subset \R^2$ with some $R>0$,  and let $m>0$ and $A>0$ be given. 
  Then there exist positive constants $T(m,A)$ and $K(m)$ such that for any
  $(u_0,v_0) \in \B (m,A)$
  the corresponding solution $(u,v)$ of \eqref{0} blows up at the finite time $\tm \le T(m,A)$.
\end{cor}

Actually, we can even give a more detailed description of nonlinear functions $\phi$ and $\psi$ yielding finite-time blowup. 
It is based on \cite[Corollary 5.2(i)]{win_mmas}.
\begin{cor}\label{cor1.35}
  If there exist $C>0$ and $s_0 >1$ such that 
  \begin{equation}\label{raz}
  \frac{\psi(s)}{\phi(s)}\geq C s\log s\;\;\mbox{for any}\;\; s>s_0>1,
  \end{equation}
  then \eqref{G1} and \eqref{H1} are satisfied. 
  Consequently, the finite-time blowup claim of Theorem~\ref{theo1} holds if \eqref{psi} and \eqref{raz} are satisfied. 
\end{cor}

Next we introduce a theorem stating the essentiality of assumption (\ref{psi}) for the dichotomy finite-time blowup 
against the boundedness of solutions. Namely we construct solutions which blow up in infinite time, once assumption (\ref{psi}) is not prescribed, although we are in the range of parameters which suggests finite-time blowup. This theorem is an extension of  \cite[Theorem 1.6]{ciesti}. Actually, for the purpose of the present paper we just need its two-dimensional part. However, since the proof is the same in higher dimensions, we give the version of the theorem which is valid for any space dimension larger than one. Let us recall that according to \eqref{0.1} the function $\beta$ is connected to $\psi$ by the identity $\psi(s) = s \beta (s)$, $s \ge 0$.

\begin{theo}\label{theo1.4}
Let $\Omega = B_R \subset \R^n$ with some $n\ge 2$ and $R>0$. Assume that there are $D_1 >0$ and 
$\gamma_1 >n$ such that for any $s\ge0$
\begin{equation}\label{balance}
  \frac{\beta^2 (s)}{\phi(s)} \le D_1 (1+s)^{-\gamma_1} 
\end{equation}
is satisfied. Moreover, suppose that there exist constants $C_1, C_2>0$ and $l_1,l_2 \in \mathbb{R}$ such that $\phi$ and $\beta$ satisfy
\begin{equation}\label{phibeta}
\phi(s)\geq C_1 (1+s)^{l_1}\;\;\;\mbox{and}\;\;\; \beta(s)\leq C_2 (1+s)^{l_2} \;\;\;\mbox{for any}\;s\ge0.
\end{equation}
Then there exists a global-in-time radially symmetric solution $(u,v)$ to \eqref{0}.

Furthermore, if additionally \eqref{G1} and \eqref{H1} are fulfilled and $n=2$, then there is a global-in-time radially symmetric solution $(u,v)$ to \eqref{0} which blows up in infinite time with respect to the norm in $L^\infty(\Omega)$.        
\end{theo}

The goal of the following remark is to show that without assuming (\ref{psi}) we still can choose such $\phi$ and $\beta$ that
\eqref{G1}, \eqref{H1} and the assumptions of Theorem~\ref{theo1.4} are satisfied at the same time. Hence, we show that (\ref{psi}) really restricts possible asymptotic behaviour of solutions to \eqref{0} to the dichotomy finite-time blowup or boundedness.
\begin{rem}\label{Uwaga}
Let $n=2$ and choose $\phi(s)=(1+s)^{-\gamma_1-2\gamma_2}$ and $\beta(s)=(1+s)^{-\gamma_1-\gamma_2}$ with some $\gamma_1 >2$ and $\gamma_2 \in (0,1)$. Then \eqref{balance} as well as \eqref{phibeta} and at the same time \eqref{G1} and \eqref{H1} are satisfied. 
\end{rem}

\mysection{Strategy of the proof of finite-time blowup}\label{zweite}

The main idea we use is a recent method introduced by M.~Winkler in \cite{win_bu}. He used it in order to treat the case of the semilinear Keller-Segel system in dimensions $n\geq 3$. We extended his method to the quasilinear system \eqref{0} in \cite{ciesti}. In the present paper we will frequently refer to the results of the latter paper. This is due to the fact that we will need 
to present some lemmata in a very precise way emphasizing the dependence of the estimates on some constants. However, the original idea and the basic estimates appeared for the first time in \cite{win_bu}.

Here we extend the method to be able to treat also a two-dimensional case. Actually, we only need to improve a single lemma. All the other parts of the proof remain the same as in \cite{ciesti}. Let us now describe the steps of the proof more precisely, often referring the reader either to \cite{ciesti} or to \cite{win_bu}.

The blowup is found as a consequence of the blowup of the Liapunov functional ${\cal F}$ associated to \eqref{0}. Namely we will show an inequality of the form 
\[
\frac{d}{dt}\left(-{\cal F} (u(t), v(t)) \right)\geq c\left(-{\cal F}^{\frac{1}{\theta}} (u(t), v(t))-1\right)
\]  
for $t>0$ with some $\theta \in (0,1)$ and $c>0$. This inequality causes blowup of ${\cal F}$ in finite time provided the initial value ${\cal F} (u_0, v_0)$ of the Liapunov functional is small enough. Hence, once we provide initial data satisfying the latter we make sure that $u$ blows up.

In order to be more precise we introduce the following notation. We fix $m>0$, $M>0$, $B>0$, and $\kappa=2$ and 
assume that
\be{m}
	\io u = m \qquad \mbox{and} \qquad \io v \le M
\ee
and 
\be{B}
	v(x) \le B|x|^{-\kappa} \qquad \mbox{for all } x\in\Omega
\ee
are fulfilled. Furthermore, we define the space
\bea{S}
	\set &:=& \bigg\{ (u,v) \in C^1(\bar\Omega) \times C^2(\bar\Omega) \ \bigg| \
	\mbox{$u$ and $v$ are positive and radially} \nn\\
	& & \hspace*{5mm}
	\mbox{symmetric satisfying $\frac{\partial v}{\partial\nu}=0$ on $\pO$, \eqref{m}, and \eqref{B}} \bigg\}.
\eea
Next we define
\be{f}
	f:=-\Delta v + v - u
\ee
and
\be{g}
	g:= \left( \frac{\phi (u)}{\sqrt{\psi(u)}}\nabla u -\sqrt{\psi (u)}\nabla v \right) \cdot \frac{x}{|x|}, \qquad x \neq 0,
\ee
for $(u,v)\in \set$.

The goal of this section is to prove that the inequality
\begin{equation}\label{4.1}
  \frac{\F(u,v)}{\D^\theta(u,v)+1} \ge - C(m,M,B) \qquad \mbox{for all } (u,v) \in \set
\end{equation} 
holds with some constants $\theta \in (0,1)$ and $C(m,M,B)>0$ . We will give the exact dependence of $C$ 
on $M$ and $B$. 

The main ingredient of the proof of \eqref{4.1} is the following estimate of $\io uv$.
\begin{lem}\label{lem3.1}
  Let \eqref{H1} and \eqref{psi} be fulfilled. Then there are $C(m)>0$ and $\theta:=\frac{8}{9}$ 
  such that all $(u,v)\in\set$ satisfy
  \begin{eqnarray}\label{3.1.1}
	\io uv &\le&  C(m) \cdot \left( 1+M^2 + B^{\frac{4}{3}} \right) \cdot 
	\Bigg( \Big\|\Delta v-v+u\Big\|_{L^2(\Omega)}^{2\theta} \nn \\
	& & + \left\|\frac{\phi (u)}{\sqrt{\psi(u)}}\nabla u -\sqrt{\psi (u)}\nabla v\right\|_{L^2(\Omega)} +1 \Bigg).
  \end{eqnarray}
\end{lem} 

Next we state two lemmata that correspond to \cite[Lemma~3.2]{ciesti} and \cite[Lemma~3.3]{ciesti}. We omit their proofs since they are exactly the same as in \cite{ciesti}, one just needs to fix $n=2$ and $\kappa=2$.

\begin{lem}\label{lem3.2}
  For any $\eps \in (0,1)$ there exists $C(\eps)>0$ such that for all $(u,v) \in \set$ 
  \be{3.2.1}
	\io uv \le (1+ \eps) \io |\nabla v|^2 + C(\eps) \cdot \left( 1+ M^2 \right) \cdot \left( 
	\Big\| \Delta v - v + u \Big\|_{L^2(\Omega)}^{\frac{4}{3}} +1 \right)
  \ee
  is fulfilled.
\end{lem}

\begin{lem}\label{lem3.3}
  For any $r_0\in (0,R)$ and $\eps \in (0,1)$, there exists a constant 
  $C(\eps,m)>0$ such that all $(u,v)\in\set$ satisfy
  \bea{3.3.1}
	\int_{\Omega \setminus B_{r_0}} |\nabla v|^2 
	&\le& \eps \io uv + \eps \io |\nabla v|^2 
    + C(\eps,m) \cdot \left( 1+ M^{\frac{4}{3}} + B^{\frac{4}{3}} \right) \cdot \Bigg\{
	r_0^{-8} \nn \\
	& & + \Big\|\Delta v-v+u\Big\|_{L^2(\Omega)}^\frac{4}{3} \Bigg\}.
  \eea
\end{lem}

Next we introduce Lemma~\ref{lem3.4} which is a main difference between \cite{ciesti} and the present paper. Both 
the statement and the proof are different to \cite[Lemma~3.4]{ciesti}. The consequence being a slightly modified, with respect to \cite{ciesti}, continuation of the proof of Lemma~\ref{lem3.1} and in turn the proof of Theorem~\ref{theo1}. 

\begin{lem}\label{lem3.4}
  Assume that \eqref{H1} and \eqref{psi} are satisfied. Then there exists  
  $C(m)>0$ such that for all $r_0\in (0,R)$ and $(u,v)\in\set$ 
  \begin{eqnarray}\label{3.4.1}
	\int_{B_{r_0}} |\nabla v|^2
	&\le&  C(m) \cdot \Bigg\{
	r_0 \cdot \Big\|\Delta v-v+u \Big\|_{L^2(\Omega)}^2 \nn \\
	& & +\left\| \frac{\phi (u)}{\sqrt{\psi(u)}}\nabla u -\sqrt{\psi (u)}\nabla v \right\|_{L^2(\Omega)} 
	+\|v\|_{L^2(\Omega)}^2
	+1 \Bigg\}
  \end{eqnarray}
  is fulfilled.
\end{lem}

The proof of Lemma~\ref{lem3.4} is the main ingredient of the next section. Another ingredient is the forthcoming 
lemma which corresponds to \cite[Lemma~3.5]{ciesti} and depends on the formulation of Lemma~\ref{lem3.4}. As its form and proof differ slightly from \cite[Lemma~3.5]{ciesti}, we shall also give its proof in the next section for the reader's convenience.   

\begin{lem}\label{lem3.5}
  Suppose that \eqref{H1} and \eqref{psi} are fulfilled and let $\theta = \frac{8}{9}$. Then for any $\eps \in (0, \frac{1}{2})$ there exists $C(\eps,m)>0$ such that 
  \begin{eqnarray}\label{3.5.1}
	\io |\nabla v|^2 
	&\le&  C(\eps,m) \cdot \left( 1+ M^2 
    + B^{\frac{4}{3}} \right) \cdot 
	\bigg( \Big\|\Delta v-v+u\Big\|_{L^2(\Omega)}^{2\theta} \nn \\ & &
	+ \Big\|\frac{\phi (u)}{\sqrt{\psi(u)}}\nabla u -\sqrt{\psi (u)}\nabla v\Big\|_{L^2(\Omega)} +1 \bigg) 
	+ \frac{\eps}{1-2\eps} \io uv
  \end{eqnarray}
  is fulfilled for all $(u,v) \in \set$.
\end{lem}

In view of Lemma~\ref{lem3.5} we are able to prove Lemma~\ref{lem3.1}.

\noindent{\bf Proof of Lemma~\ref{lem3.1}.} We fix some $\eps \in (0, \frac{1}{2})$ and apply Lemma~\ref{lem3.2} to deduce that 
$$\io uv \le (1+ \eps) \io |\nabla v|^2 + c_1 \cdot \left( 
	\| f \|_{L^2(\Omega)}^{\frac{4}{3}} +1 \right)$$
is satisfied with $c_1 = C_1 \cdot (1+M^2)>0$.	
Furthermore, Lemma~\ref{lem3.5} implies the existence of $c_2 = C_2(m) \cdot \big( 1+ M^2 
    + B^{\frac{4}{3}} \big)>0$ such that
\[
\io uv \le   \frac{\eps(1+ \eps)}{1-2\eps} \io uv 
   + c_2(1+ \eps) \cdot 
	\left(\|f\|_{L^2(\Omega)}^{2\theta} 
	+ \|g\|_{L^2(\Omega)} +1 \right) + c_1 \cdot \left( 
	\| f \|_{L^2(\Omega)}^{\frac{4}{3}} +1 \right)
\]
which yields 
\[
\io uv \le c_3 \cdot 
	\left(\|f\|_{L^2(\Omega)}^{2\theta} + \| f \|_{L^2(\Omega)}^{\frac{4}{3}}
	+ \|g\|_{L^2(\Omega)} +1 \right) 
\]
with some $c_3 = C_3(m) \cdot \big( 1+ M^2 
    + B^{\frac{4}{3}} \big)>0$. In view of  $\frac{4}{3}<2\theta=\frac{16}{9}$, 
a further application of the Young inequality implies \eqref{3.1.1}.
\qed
Finally, we show that the Liapunov functional $\F$ can be estimated according to \eqref{4.1}. 
\begin{theo}\label{theo3.6}
  Assume that \eqref{H1} and \eqref{psi} are satisfied and let $\theta=\frac{8}{9}$. Then there exists $C(m)>0$
  such that 
  \be{3.6.1}
	\F(u,v) \ge -C(m) \cdot \left( 1+M^2 + B^{\frac{4}{3}} \right) \cdot \Big( \D^\theta (u,v) +1 \Big)
  \ee
  is fulfilled for all $(u,v) \in \set$, where $\F$ and $\D$ are given in \eqref{F} and \eqref{D}, 
  respectively.
\end{theo}
\proof
  In view of \eqref{f}, \eqref{g}, and $\theta > \frac{1}{2}$, an application of Young's inequality
  to \eqref{3.1.1} implies the existence of $c_1=C_1(m) \cdot \big( 1+M^2 + B^{\frac{4}{3}} \big)>0$ such that
  \bas
	\io uv \le c_1 \left( \Big(\|f\|_{L^2(\Omega)}^2 + \|g\|_{L^2(\Omega)}^2 \Big)^\theta +1 \right).
  \eas
  As moreover \eqref{0.1} and \eqref{GH} imply that $G$ is nonnegative, we deduce that
  \bas
	\F(u,v)
	&=& \frac{1}{2} \io |\nabla v|^2 + \frac{1}{2} \io v^2 - \io uv + \io G(u) \\
	&\ge& - c_1 \cdot \left( \left(\|f\|_{L^2(\Omega)}^2 + \|g\|_{L^2(\Omega)}^2\right)^\theta +1 \right).
  \eas
  Since \eqref{D}, \eqref{f}, and \eqref{g} imply $\D(u,v)=\|f\|_{L^2(\Omega)}^2 + \|g\|_{L^2(\Omega)}^2$,
  the claim is proved.
\qed
Now we are in a position to prove the finite-time blowup of solutions to \eqref{0}.

{\bf Proof of Theorem~\ref{theo1}.} Since the proof of \cite[Corollary~3.3]{win_bu} is based on estimates 
coming only from the second equation of \eqref{0} and it is not changed for dimension $n=2$, the corollary remains true with $\kappa=2$. Consequently, 
we know that $v(t)$ satisfies \eqref{B} for all $t \in (0, \tm)$ with 
\[
B\leq C \left(\|u_0\|_{L^1(\Omega)}+\|v_0\|_{L^1(\Omega)}+\|\nabla v_0\|_{L^2(\Omega)}\right).
\]    
Next the proof of Theorem~\ref{theo1} splits into two parts. The part of Theorem~\ref{theo1} concerning finite-time blowup of solutions 
provided they start from initial data belonging to ${\cal B}(m,A)$ follows exactly the lines of \cite[Lemma~5.2]{win_bu}. The exact dependence of ${\cal F} (u_0,v_0)$ on $A$ can be shown like in \cite[Lemma~4.1]{ciesti}. Furthermore, given an arbitrary $m>0$, \cite[Lemma 4.1]{win_mmas} guarantees the existence of $A>0$ such that the set ${\cal B}(m,A)$ is nonempty. Indeed, 
choosing the functions $(u_\eta,v_\eta)$, $\eta >0$, which are defined there, we see that for $\eta$ small enough    
$\F(u_\eta,v_\eta) \le -K(m) \cdot (1+A_\eta^2)$ is satisfied with $A_\eta=\| v_\eta\|_{W^{1,2}(\Omega)}$ and conclude that
${\cal B}(m,A_\eta) \neq \emptyset$. 
\qed

\mysection{Main estimates}\label{dritte}
The present section is devoted to proving Lemma \ref{lem3.4} which is the main contribution of our paper with respect to the estimates showing finite-time blowup. Moreover, 
we give the proof of Lemma \ref{lem3.5} which is just a slight modification of \cite[Lemma 3.5]{ciesti} due to a different formulation of the preceding lemma. 

\noindent
{\bf Proof of Lemma \ref{lem3.4}} Since $u$ and $v$ are radially symmetric, \eqref{f} and \eqref{g} 
imply 
\begin{equation}\label{pier}
(rv_r)_r = -ru - rf + rv
\end{equation}
and 
\begin{equation}\label{pierd}
v_r = \frac{\phi(u)}{\psi(u)} u_r- \frac{g}{\sqrt{\psi(u)}}\;.
\end{equation} 
Next, by multiplying (\ref{pier}) by $rv_r$, using (\ref{pierd}) and applying 
Young's inequality, for any $\delta >0$ we deduce that 
\begin{eqnarray}\label{3.4.2}
 \frac{1}{2} \left( (r v_r)^2 \right)_r &=& -r^2uv_r - r^2 fv_r + r^2 vv_r \nn \\
 &\le& -r^2 \frac{u \phi(u)}{\psi(u)} u_r + r^2 \frac{u}{\sqrt{\psi(u)}} g 
 + \frac{\delta}{2} (rv_r)^2 + \frac{1}{2 \delta} r^2 f^2 \nn \\
 & & + \frac{1}{2} r^2 (v^2)_r \qquad \mbox{for all } r \in (0,R).
\end{eqnarray}
Moreover, we define $y(r) := (r v_r)^2$, $r \in [0,R]$, and obtain 
$$y_r \le -2r^2 \frac{u \phi(u)}{\psi(u)} u_r + 2r^2 \frac{u}{\sqrt{\psi(u)}} g 
 + \delta y + \frac{1}{\delta} r^2 f^2 + r^2 (v^2)_r, \qquad r \in (0,R),$$ 
as well as $y(0) =0$ due to the regularity of $v$. Thus, an integration yields
\begin{eqnarray}\label{3.4.3}
 r^2 v_r^2(r) = y(r) &\le&
	-2\int_0^r e^{\delta (r-\rho)} \rho^2 \frac{u(\rho) \phi(u(\rho))}{\psi(u(\rho))} u_r (\rho) 
        \; d\rho \nn \\ & &
	+2\int_0^r e^{\delta (r-\rho)} \rho^2 \frac{u(\rho)}{\sqrt{\psi(u(\rho))}} g (\rho) \; d\rho \nn\\
	& & +\frac{1}{\delta} \int_0^r e^{\delta (r-\rho)} \rho^2 f^2(\rho) d\rho
	+\int_0^r e^{\delta (r-\rho)} \rho^2 (v^2)_r(\rho) d\rho
\end{eqnarray}
for all $r \in (0,R)$.
Integrating by parts and using the nonnegativity of $H$, we estimate
\begin{eqnarray}\label{3.4.4}
  && \hspace*{-20mm} 
  -2\int_0^r e^{\delta (r-\rho)} \rho^2 \frac{u (\rho) \phi(u (\rho))}{\psi(u (\rho))} u_r (\rho) \; d\rho 
  \nn \\
  &=& 4 \int_0^r e^{\delta (r-\rho)} \rho H(u(\rho))  \; d\rho \nn \\
  & & -2 \delta \int_0^r e^{\delta (r-\rho)} \rho^2 H(u(\rho)) \; d\rho
  -2r^2 H(u(r)) \nn \\
  &\le& 4 e^{\delta R} \int_0^r \rho H(u(\rho))  \; d\rho, \qquad r \in (0,R).
\end{eqnarray}
Next, by the 
Cauchy-Schwarz inequality and \eqref{psi} we deduce that 
\begin{eqnarray}\label{3.4.5}
 && \hspace*{-20mm} 
 2\int_0^r e^{\delta (r-\rho)} \rho^2 \frac{u(\rho)}{\sqrt{\psi(u(\rho))}} g (\rho) \; d\rho \nn \\
 &\le& 2 \left( \int_0^R \rho \frac{u^2(\rho)}{\psi(u(\rho))} \; d\rho \right)^\frac{1}{2} \cdot 
	\left( \int_0^r e^{2\delta(r-\rho)} \cdot \rho^3 g^2(\rho) \; d\rho \right)^\frac{1}{2} \nn \\
 &\le& 2 \left( \frac{1}{c_0} \int_0^R \rho u(\rho) \; d\rho \right)^\frac{1}{2} \cdot 
	\left( e^{2\delta R} r^2 \int_0^R \rho g^2(\rho) \; d\rho \right)^\frac{1}{2} \nn \\ 
 &\le& \frac{re^{\delta R}}{\pi \sqrt{c_0}}  \sqrt{m}\| g \|_{L^2(\Omega)},
 \qquad r \in (0,R).
\end{eqnarray}
Similarly, we estimate the third term on the right-hand side of \eqref{3.4.3} according to 
\begin{eqnarray}\label{3.4.6}
  \frac{1}{\delta} \int_0^r e^{\delta (r-\rho)} \rho^2 f^2(\rho) \; d\rho
	&\le& \frac{re^{\delta R}}{\delta} \int_0^R \rho f^2(\rho) \; d\rho \nn\\
	&=& \frac{re^{\delta R}}{2\pi\delta}  \|f\|_{L^2(\Omega)}^2
	\qquad \mbox{for all } r\in (0,R).
\end{eqnarray}
Now we fix $\delta \in (0, \frac{2}{R})$. Hence, $2 \rho \ge \delta \rho^2$ for all $\rho \in (0,R)$,
and an integration by parts yields 
\begin{equation}\label{3.4.7}
 \int_0^r e^{\delta(r-\rho)} \rho^2 (v^2)_r(\rho) \; d\rho
	\le r^2 v^2(r)
	\qquad \mbox{for all } r\in (0,R).
\end{equation}
Thus, \eqref{3.4.3}-\eqref{3.4.7} imply that there is a constant $c_1(m) >0$ such that
\begin{eqnarray*}
  r^2 v_r^2(r) &\le& 4 e^{\delta R} \int_0^r \rho H(u(\rho))  \; d\rho \\
  & & + \frac{c_1(m)}{2\pi} r \|g\|_{L^2(\Omega)} 
   + \frac{c_1(m)}{2\pi} r \|f\|_{L^2(\Omega)}^2
	+ r^2 v^2(r), \qquad r \in (0,R).
\end{eqnarray*}
Multiplying this inequality by $2\pi r^{-1}$ and integrating over $r \in (0,r_0)$, we conclude that
\begin{eqnarray}\label{3.4.8}
 \int_{B_{r_0}} | \nabla v |^2 &=& 2\pi \int_0^{r_0} r v_r^2(r) \; dr \nn \\
 &\le& 8\pi e^{\delta R}  \int_0^{r_0} r^{-1} \int_0^r \rho H(u(\rho))  \; d\rho \, dr \nn \\
 & & +c_1 (m)R \|g\|_{L^2(\Omega)}+ c_1(m) r_0 \|f\|_{L^2(\Omega)}^2 + \|v\|_{L^2(\Omega)}^2. 
\end{eqnarray}
In order to estimate the first term on the right-hand side, we adapt an idea used in \cite[Lemma~3.3]{win_mmas} and first claim that
\begin{equation}\label{33}
  H(s) \ln (H(s)) \le c_2 (1+s), \qquad s \ge 0,
\end{equation}
is fulfilled with some positive constant $c_2$. 
In view of \eqref{H1} there exists $c_3 >0$ such that
\begin{eqnarray*}
  H(s) \ln (H(s)) &\le& b \frac{s}{\ln s} \cdot \ln \left( \frac{bs}{\ln s} \right)
  = b \frac{s}{\ln s} \cdot \left( \ln s + \ln \left( \frac{b}{\ln s} \right) \right) 
  \le c_3 (1+s), \qquad s \ge s_0,  
\end{eqnarray*}
which implies \eqref{33} due to \eqref{0.1} and the definition of $H$ in \eqref{GH}.
Let us further remind that Young's inequality
\begin{equation}\label{22}
 AB \le \frac{1}{e} e^A + B \ln B
\end{equation}
holds for all $A, B >0$. Applying now Fubini's theorem to the first term on the right-hand side of (\ref{3.4.8})
and using (\ref{22}) and (\ref{33}), we obtain 
\begin{eqnarray*}
  && \hspace*{-20mm}
   8\pi e^{\delta R}  \int_0^{r_0} r^{-1} \int_0^r \rho H(u(\rho))  \; d\rho \, dr \\
  &=& 8 \pi e^{\delta R} \int_0^{r_0} \left( \int_{\rho}^{r_0} r^{-1} \; dr \right) 
    \rho H(u(\rho))  \; d\rho \\
  &=& 8 \pi e^{\delta R} \int_0^{r_0} \ln \left( \frac{r_0}{\rho} \right)
    \rho H(u(\rho))  \; d\rho \\
  &\le& 8 \pi e^{\delta R} \int_0^{r_0} \left(  \frac{1}{e} \cdot \frac{r_0}{\rho} \cdot \rho  
    + \rho H(u(\rho)) \ln (H(u(\rho))) \right)  \; d\rho \\
  &\le& 8 \pi e^{\delta R} \int_0^{r_0} \left(  \frac{r_0}{e} 
    + \rho c_2 (1+ u(\rho)) \right) \; d\rho 
  \le c_4 r_0^2 + c_4 \int_{B_{r_0}} u
  \le c_4 R^2 + c_4 m
\end{eqnarray*}
with some $c_4 >0$. In the light of (\ref{3.4.8}) the lemma is proved.
\qed
{\bf  Proof of Lemma \ref{lem3.5}.} Let us fix $r_0:= \min \{\frac{R}{2},\|f\|_{L^2(\Omega)}^{-\frac{2}{9}}\}$.
By applying Lemma~\ref{lem3.3} we obtain
$c_1= C_1 (\eps, m) \cdot \big( 1+ M^{\frac{4}{3}} + B^{\frac{4}{3}} \big)>0$ such that
\begin{equation}\label{3.5.2}
 \int_{\Omega \setminus B_{r_0}} |\nabla v|^2
	\le \eps \io uv + \eps \io |\nabla v|^2
	+ c_1 \cdot \Big( r_0^{-8} + \|f\|_{L^2(\Omega)}^\frac{4}{3} \Big).
\end{equation}
Next, by Lemma~\ref{lem3.4} there exists $c_2 = c_2(m)>0$ such that
\begin{equation}\label{3.5.3}
 \int_{B_{r_0}} |\nabla v|^2  \le c_2 \cdot \Big( r_0 \|f\|_{L^2(\Omega)}^2 + \|g\|_{L^2(\Omega)} 
   + \|v\|_{L^2(\Omega)}^2 + 1 \Big).
\end{equation}
As moreover 
$$c_2 \|v\|_{L^2(\Omega)}^2 \le \eps \io |\nabla v|^2 + c_3$$
is fulfilled by \eqref{m} and \cite[Lemma~2.2]{win_bu} with some $c_3 = C_3 (\eps,m) \cdot M^2 >0$,
by adding (\ref{3.5.2}) and (\ref{3.5.3}) we deduce that
\begin{eqnarray}\label{3.5.5}
 (1-2\eps) \io |\nabla v|^2  &\le& \eps \io uv + c_2 (\|g\|_{L^2(\Omega)} +1) + c_3  +I,
\end{eqnarray}
where  
$$I := c_1 r_0^{-8} + 
  c_1 \|f\|_{L^2(\Omega)}^\frac{4}{3}  + c_2 r_0 \|f\|_{L^2(\Omega)}^2.$$
In case of $\|f\|_{L^2(\Omega)} \le (\frac{2}{R})^\frac{9}{2}$, we have $r_0=\frac{R}{2}$ and hence $I$ can be estimated by a constant. This proves \eqref{3.5.1} in this case.

Furthermore, in the case $\|f\|_{L^2(\Omega)} > (\frac{2}{R})^\frac{9}{2}$ we have 
$r_0=\|f\|_{L^2(\Omega)}^{-\frac{2}{9}}$ which implies  
$$ I \le c_1 \|f\|_{L^2(\Omega)}^\frac{16}{9}
	+c_1 \|f\|_{L^2(\Omega)}^\frac{4}{3}
	+ c_2 \|f\|_{L^2(\Omega)}^{\frac{16}{9}} 
	= (c_1+c_2) \|f\|_{L^2(\Omega)}^\frac{16}{9} 
	+c_1 \|f\|_{L^2(\Omega)}^\frac{4}{3}.$$
Since $\frac{4}{3}< \frac{16}{9}=2\theta$ we are in the position to use Young's inequality once more to obtain 
$$I \le (2c_1+c_2) \|f\|_{L^2(\Omega)}^{2\theta} +  c_1,$$
which inserted into \eqref{3.5.5} proves \eqref{3.5.1} in the case 
  $\|f\|_{L^2(\Omega)} > (\frac{2}{R})^\frac{9}{2}$ and thereby completes the proof.
\qed

\mysection{Infinite-time blowup}\label{vierte}
This section is devoted to the proof of Theorem \ref{theo1.4}. To this end we first prove the following lemma which generalizes \cite[Lemma~5.1]{ciesti}.
\begin{lem}\label{lemat}
Let $\Omega\subset \R^n$ be a bounded domain with some $n\ge 2$. Moreover, assume that \eqref{balance} holds. 
Then there exists $p>n$ such that for any solution $(u,v)$ to \eqref{0} and any $T \in (0, \infty)$ with $T \le \tm$ there is $C>0$ such that $u$ admits the estimate
\begin{equation}\label{wzor}
\|u(\cdot,t)\|_{L^p(\Omega)}\leq C, \qquad t \in \left( 0,T \right).
\end{equation}
\end{lem}
\proof We fix $p \in (n, \gamma_1]$, multiply the first equation of \eqref{0} by $u^{p-1}$ and the second one by
$\Delta v$ in order to obtain
\begin{equation}
\label{eins}
\frac{1}{p}\frac{d}{dt}\int_\Omega u^pdx + (p-1)\int_\Omega \phi(u)\left|\nabla u\right|^2u^{p-2}dx = (p-1)\int_\Omega u^{p-1}\beta(u)\nabla v\nabla u \; dx
\end{equation}
and 
\begin{equation}
\label{zwei}
\frac{1}{2}\frac{d}{dt}\int_\Omega |\nabla v|^2dx + \frac{1}{2}\int_\Omega\left|\Delta v\right|^{2}dx+ \int_\Omega|\nabla v|^2dx\leq \frac{1}{2}\int_\Omega u^2dx.
\end{equation}
Writing
\[
u^{p-1}\beta(u)=u^{\frac{p-2}{2}} \sqrt{\phi(u)} u^{\frac{p}{2}}\frac{\beta(u)}{\sqrt{\phi(u)}},
\]
we deduce from (\ref{eins}) that
\begin{equation}\label{drei}
\frac{1}{p}\frac{d}{dt}\int_\Omega u^pdx + \frac{p-1}{2}\int_\Omega \phi(u)\left|\nabla u\right|^2u^{p-2}dx\leq
C \int_\Omega u^p\frac{\beta^2(u)}{\phi(u)}|\nabla v|^2dx.
\end{equation}
Next adding (\ref{drei}) and (\ref{zwei}), applying (\ref{balance}) and using $p \le \gamma_1$, we conclude that
\begin{equation}\label{vier}
\frac{d}{dt}\left(\int_\Omega u^pdx+\int_\Omega|\nabla v|^2dx\right)\leq C\left(\int_\Omega u^pdx\right)^{\frac{2}{p}}+C\int_\Omega|\nabla v|^2\leq C\left(\int_\Omega u^pdx+\int_\Omega|\nabla v|^2dx +1\right).
\end{equation}
Now Gr\"{o}nwall's lemma implies the claimed estimate of $\|u\|_{L^p(\Omega)}$. 
\qed  
Now we can prove the blowup in infinite time by a suitable combination of known results.

{\bf Proof of Theorem \ref{theo1.4}.} Due to (\ref{wzor}) and the classical regularity theory of parabolic equations
applied to the second equation of \eqref{0}, see \cite[Lemma 4.1]{horstmann_winkler} for example, one obtains an estimate of $\nabla v$ in $L^\infty (\Omega \times (0,T))$ for any finite $T \in (0, \tm]$. Next we multiply the first equation of \eqref{0} by $u^{p-1}$, this time for any $p \in (\gamma_1,\infty)$. Proceeding as in the proof of Lemma~\ref{lemat}, we see that the right-hand side of (\ref{drei}) can be estimated by $C(\|\nabla v\|_{L^\infty(\Omega)})\io u^{p-\gamma_1}$ due to (\ref{balance}). Hence, H\"{o}lder's inequality leads to
\[
\frac{d}{dt}\int_\Omega u^pdx\leq C\left(\int_\Omega u^pdx+1\right).
\] 
Thus, in view of Lemma~\ref{lemat}, for any $p \in (1,\infty)$ the norm $\|u\|_{L^p(\Omega)}$ can be bounded according to 
\eqref{wzor}. We are now in a position to apply \cite[Lemma A.1]{taowin_jde} in order to gain an estimate of $u$ in $L^\infty
(\Omega \times (0,T))$ which shows the existence of a global solution. More precisely, keeping the notation of \cite[Lemma A.1]{taowin_jde}, we have $f:=u\beta(u)\nabla v$ and $g:=0$, while due to (\ref{phibeta})
we can choose $m=l_1 +1$. Moreover, by \eqref{phibeta} and the estimates on $u$ we just proved, we have $u \in L^\infty ((0,T); L^{p_0} (\Omega))$ and $f \in L^\infty ((0,T); L^{q_1} (\Omega))$ for any $p_0 \in (1,\infty)$ and $q_1 \in (1,\infty)$. This freedom of choosing any $p_0<\infty$ as well as any $q_1<\infty$ enables us to make sure that all the assumptions of  \cite[Lemma A.1]{taowin_jde} are satisfied.

Furthermore, if we additionally assume that \eqref{G1} and \eqref{H1} are satisfied with $n=2$, we apply 
\cite[Theorem~5.1]{win_mmas} in order to deduce that $(u,v)$ blows up in infinite time. This finishes the proof of 
Theorem~\ref{theo1.4}.
\qed

\end{document}